\documentclass[a4paper,12pt]{article}
\usepackage[utf8]{inputenc}
\usepackage{multicol}
\usepackage{a4}
\usepackage{amsthm,amsmath,amssymb,setspace}
\usepackage{pgf,tikz,pgfplots}
\usetikzlibrary{arrows}
\usepackage{mathrsfs}
\usepackage{graphicx}
\usepackage{array}
\usepackage{hyperref}
\usepackage{titlesec}
\usepackage{enumerate}
\usepackage[english]{babel}
\usepackage{longtable}
\usepackage{enumerate}
\usepackage{hyperref}
\usepackage{subfigure}
\usepackage[linesnumbered,ruled,vlined]{algorithm2e}
\usepackage{cite}
\usepackage{mathtools}
\usepackage{setspace}
\usepackage{algorithmic}
\usepackage{lineno}
\usepackage{rotating}
\usepackage{float}
\usepackage{lscape}
\usepackage{amsmath}
\SetKwInput{KwInput}{Input}
\SetKwInput{KwOutput}{Output}
\newtheorem{theorem}{Theorem}[section]

\begin{document}

\title{Antimagic Labeling of Generalized Edge Corona Graphs}
\author{\textsc{D. Nivedha, S. Devi Yamini}\\{\footnotesize Vellore Institute of Technology, Chennai, Tamil Nadu, India}\\ {\footnotesize\texttt{email: nivedha.d2020@vitstudent.ac.in, deviyamini.s@vit.ac.in}}}

\date{ }
\maketitle

\begin{abstract}
An antimagic labeling of a graph $G$ is a one-to-one  correspondence between the edge set $E(G)$ and $\lbrace 1,2,...,|E(G)|\rbrace$ in which the sum of the edge labels incident on the distinct vertices are distinct. Let  $G$,$H_1$,$H_2$,...,$H_{m-1}$, and $H_m$ be simple graphs where $|E(G)|=m$. A generalized edge corona of the graph $G$  and $(H_1,H_2,...,H_m)$ (denoted by $G\diamond (H_1,H_2,...H_m)$) is a graph obtained by taking a copy of $G,H_1,H_2,...,H_m$ and joining the end vertices of $i^{th}$ edge of $G$ to every vertex of $H_i$, $i\in\lbrace 1,2,...,m\rbrace$. In this paper, we consider $G$ as a connected graph with exactly one vertex of maximum degree 3 (excluding the spider graph with exactly one vertex of maximum degree 3 containing uneven legs) and each $H_i$, $1\leq i \leq m$ as a connected graph on at least two vertices. We provide an algorithmic approach to prove that $G$ $\diamond$ $(H_1,H_2,...H_m)$ is antimagic under certain conditions.\\
\end{abstract}
\textbf{Keywords: Graph labeling, Antimagic labeling, Generalized edge corona graphs, Spider graphs, Pan graphs}

	\section{Introduction}
	The graphs involved in this paper are simple and undirected. Most of the real life situations can be modeled as a graph problem which has attracted the researchers to work on graph theory. Graph labeling is one of the topics in graph theory which has more than 200 techniques. The concept of graph labeling was introduced by Alexander Rosa in 1967. A graph labeling is an assignment of integers to the vertices or edges (or both) subject to certain conditions. Refer \cite{23} for more details on graph labeling. Hartsfield and Ringel introduced the concept of an antimagic labeling in the year 1990.  Antimagic labeling has numerous results and still a lot of work is under process. Let $V(G)$ and $E(G)$ denote the vertex set and the edge set of the graph $G$ respectively where $|V(G)|=n$. For a graph $G=(V,E)$, an antimagic labeling is a bijection $f:E(G)\rightarrow \lbrace 1,2,...,|E(G)| \rbrace$ such that $w(u)=\sum\limits_{e\in E(u)}f(e)$ is distinct for all vertices $u\in V(G)$ where $E(u)$ denotes the set of all edges incident on the vertex $u$. A graph is antimagic if it admits an antimagic labeling. In \cite{11}, Hartsfield and Ringel (1990) proposed the following two conjectures which remain open for more than three decades.\\
 \textbf{Conjecture 1.1.}\cite{11} Every connected graph other than $K_2$ is antimagic.\\
\textbf{Conjecture 1.2.}\cite{11} Every tree other than $K_2$ is antimagic.\\
The Table \ref{tab:my_label} provides a few existing results on antimagic labeling.
\begin{table}[!htb]
    \begin{tabular}{ | m{12cm} | m{2cm}| } 
   
      \hline
      \centering \textit{Graphs} & \textit{Ref.}\\
    \hline
      - Paths\newline
      - Cycles \newline
      - Wheels \newline
      - Complete graphs
      &  \cite{11}    \\  \hline
      - Graphs with $\triangle(G)\geq n-3$ & \cite{12}\\ \hline
    - Toroidal grids\newline
- Higher dimensional toroidal grids\newline
- Cartesian product of cycle and $k$-regular graph &   \cite{13}\\ \hline
- Sequential
generalized corona graphs\newline
- Generalized snowflake graphs &  \cite{17}\\ \hline
- $n$-barbell graph  $n\geq 3$\newline
- Edge corona of bistar graph and $k$-regular graph\newline
- Edge corona of cycles &    \cite{20}\\ \hline
- Binomial trees\newline
- Fibonacci trees &   \cite{21} \\ \hline
Complete $m$-ary trees &  \cite{22}\\ \hline
Subclasses of trees &  \cite{9}\\ \hline
Caterpillars &    \cite{16}\\ \hline
Regular graphs &   \cite{1,2}\\ \hline
Biregular bipartite graphs & \cite{41}\\ \hline
Hexagonal lattice\newline
Prismatic lattice & \cite{31}\\ \hline

 \end{tabular}
    \caption{Existing results on antimagic labeling}
    \label{tab:my_label}
\end{table}
\section{Motivation and Applications }
A few classes of graphs and products of graphs were proved to be antimagic, but the antimagic labeling is yet to be explored for the edge corona product of graphs. So, this motivated us to work on the generalized edge corona graph (denoted by $G$ $\diamond$ $(H_1,H_2,...,H_m)$).\\

\indent There are techniques used in surveillance or security model for various buildings which are based on antimagic labeling of double wheel graph, centreless wheel graph, helm, and regular actinia graphs \cite{25}. A few antimagic graphs such as double wheel, helm, path, web, etc. are used in encryption techniques for the security purpose in data transfer \cite{26}.

\section{Preliminaries}
\begin{figure}[!htb]
\centering		\includegraphics[width=1.05\linewidth,height=0.2\textheight]{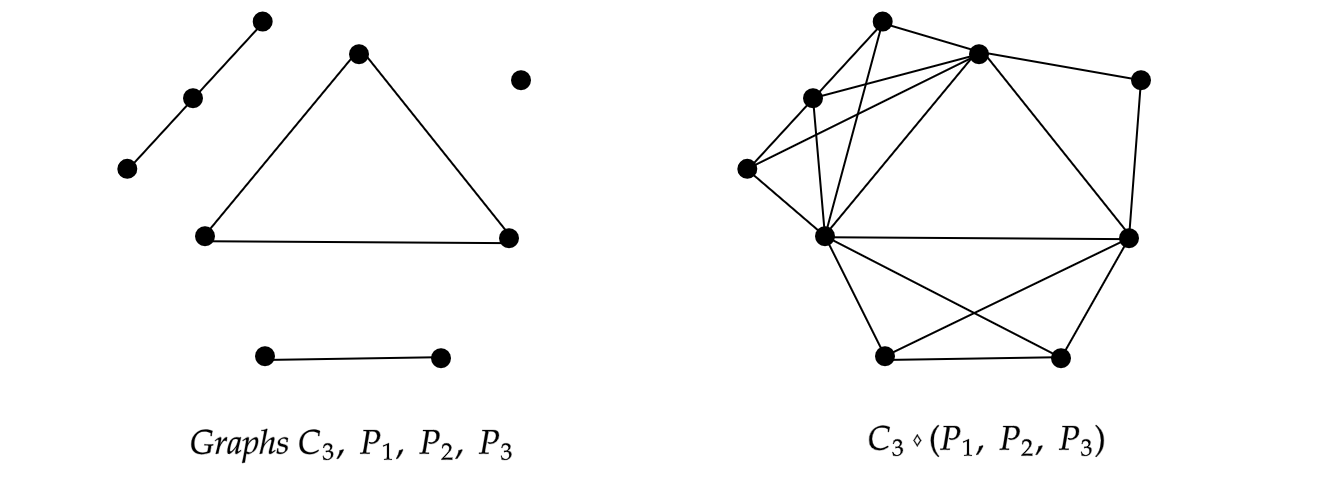}
		\caption{A generalized edge corona of graphs}\label{fig:A}
	\end{figure}

The generalized edge corona of a graph $G$ (on $m$ edges) and $(H_1,H_2,...,H_m)$ (denoted by $G\diamond (H_1,H_2,...H_m)$) is a graph obtained by taking a copy of $G,H_1,H_2,...,H_m$ and joining the end vertices of $i^{th}$ edge of $G$ to every vertex of  $H_i$, $i\in\lbrace 1,2,...,m\rbrace$ \cite{14}. An illustration of generalized edge corona of the graph $C_3$ and $(P_1,P_2,P_3)$ is given in the Figure \ref{fig:A}.
 
 The edges between the graph $G$ and $H_i$, $i\in \lbrace 1,2,...,m\rbrace$ are denoted as cross edges. The degree of the vertex $v\in V(G)$ is defined as the number of edges incident on the vertex $v$. For any vertex $v\in V(G\diamond (H_1,..H_m))$, we define
\begin{equation}
d(v) =
\begin{cases}
\text{degree of $v$ in $G$} & \text{if $v\in V(G)$}\nonumber\\
\text{degree of $v$ in $H_i$} & \text{if $v\in V(H_i)$}\nonumber
\end{cases}
\end{equation}
 and $d'(v)$ represents the degree of $v$ in $G\diamond (H_1,H_2,...H_m)$. In a similar manner, we define the maximum degrees $\triangle$, $\triangle'$ and the minimum degrees $\delta$, $\delta'$. The sum of the labels of all the edges incident on a vertex $v$ in a graph $G$ is called the vertex sum of $v$ (denoted as $w(v)$); and the sum of the labels of some edges (i.e., few edges remain unlabeled) incident on a vertex $v$ in a graph $G$ is known as the partial vertex sum of $v$ (denoted as $w'(v)$). Let PVS and VS be the acronym for partial vertex sum and the vertex sum respectively.   The join of any two graphs $G$ and $H$ is obtained by joining every vertex of $G$ to all the vertices of $H$. A pan graph is obtained by joining a vertex of a cycle graph to a singleton graph by an edge. A spider graph is a tree with a  vertex of degree at least $3$ and other vertices of degree at most $2$.

\section{ Generalized edge corona of graphs }

Let $G$ be a connected graph with exactly one vertex of maximum degree 3 (excluding the spider graph with exactly one vertex of maximum degree 3 containing uneven legs) and the other vertices of degree $2$ or $1$. Let the graphs $H_0, H_1,...,H_r$, $r\geq 3$ be connected with at least two vertices. Note that the graph $G$ can be classified into two types as follows:\\

\begin{figure}[!htb]
		\centering
		\includegraphics[width=.5\linewidth]{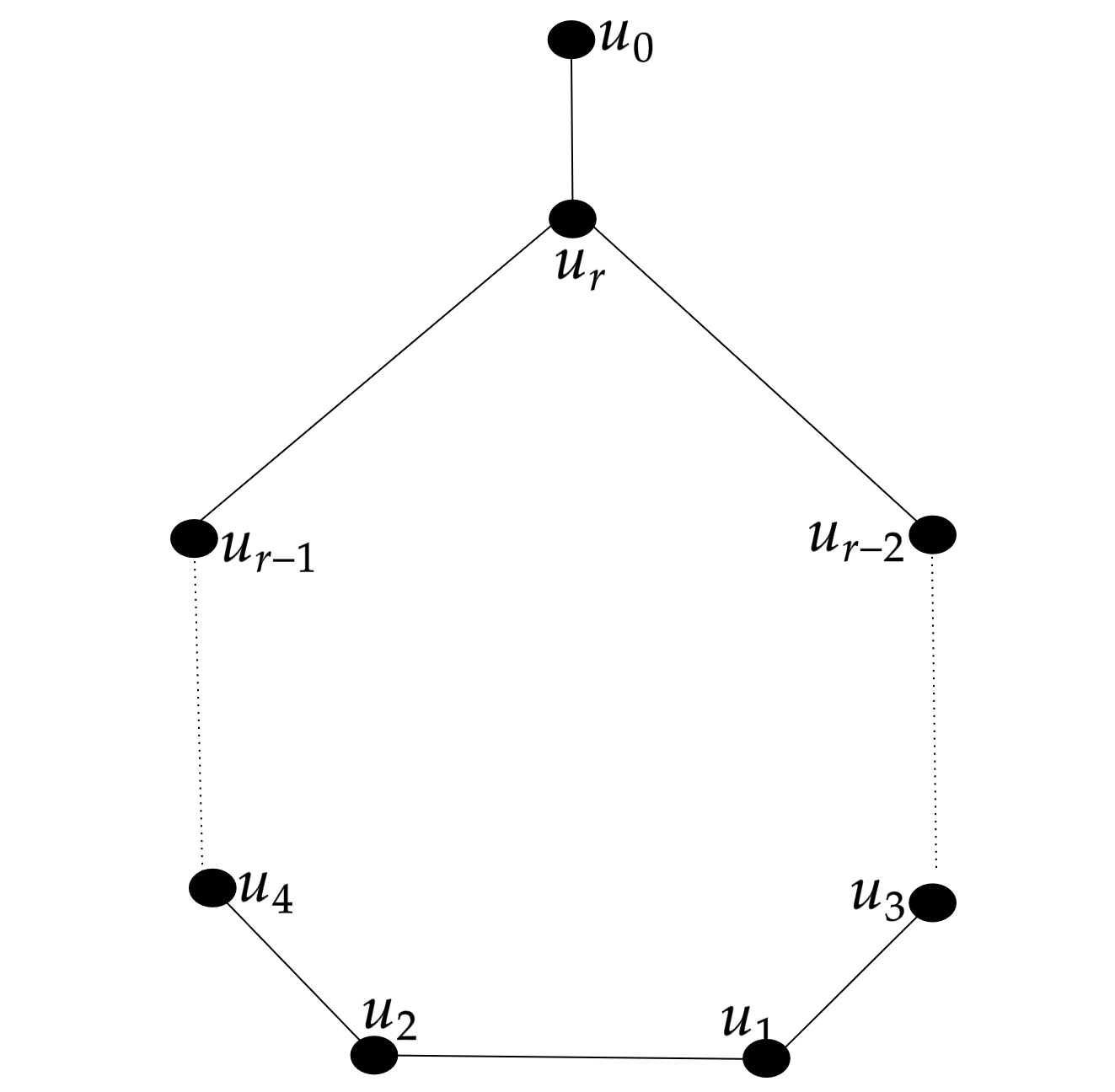}
		\caption{Type I graph $G_1$}\label{fig:1}
	\end{figure}

\textbf{Type I:} A pan graph $G_1$ with $V(G_1)=\lbrace u_0,u_1,u_2,...,u_r\rbrace$ where $u_0$ is the pendant vertex and $u_i$, $1\leq i\leq r$ ordered as in the Figure \ref{fig:1} (the dotted line from $u_4$ to $u_{r-1}$ represent a path $u_4-u_6-u_8-...-u_{r-1}$ and the dotted line from $u_3$ to $u_{r-2}$ represent a path $u_3-u_5-u_7-...-u_{r-2}$).
\noindent\\

\begin{figure}[!htb]
		\centering
		\includegraphics[width=.5\linewidth]{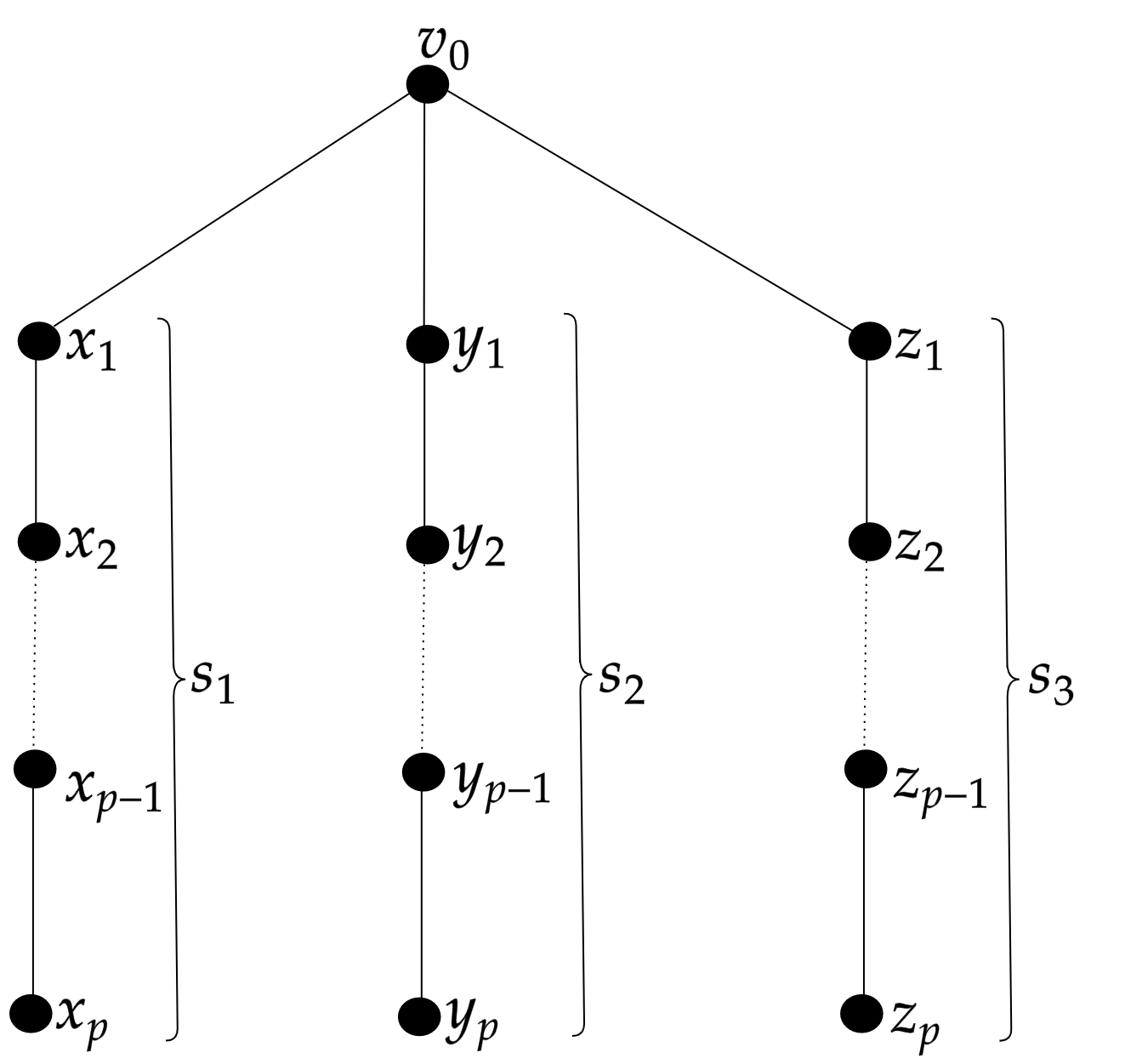}
		\caption{Type II graph $G_2$}\label{fig:2}
	\end{figure}

\textbf{Type II:} A spider graph $G_2$ with $v_0$ as the vertex of degree $3$. Let $s_i$, $1\leq i \leq 3$ be the leg of the spider where each leg represents a path on $p\geq 1$ vertices as in the Figure \ref{fig:2} (the dotted line from $x_{p-1}$ to $x_2$ represent a path $x_{p-1}-x_{p-2}-...-x_{2}$, the dotted line from $y_{p-1}$ to $y_2$ represent a path $y_{p-1}-y_{p-2}-...-y_{2}$, and the dotted line from $z_{p-1}$ to $z_2$ represent a path $z_{p-1}-z_{p-2}-...-z_{2}$). We omit the spider graph with exactly one vertex of maximum degree 3 containing uneven legs from Type II.

	\hfill
 
\subsection{$G_{1}$ $\diamond$ $(H_{0},H_{1},...,H_{r})$ is antimagic}
Note that $G_1$ is a Type I graph. Without loss of generality, assume that $|V(G_1)|=|E(G_1)|=r+1$, $r\geq 3$ and $u_0u_r\in E(G_1)$. Recall that the graph $H_i$, $i\in \lbrace 0,1,2,...,r\rbrace$ is connected with at least two vertices. Let $V(H_i)=\lbrace v^i_1,...,v^i_{n_i}\rbrace$ ($|V(H_i)|=n_i$) and $|E(H_i)|=q_i, \hspace{0.1cm} \forall\hspace{0.1cm} i\in \lbrace 0,1,...,r\rbrace$. Arrange the graphs $H_i$ satisfying the following condition: $|V(H_i)|\leq |V(H_{i+1})|$, $\forall$ $i\in \lbrace0,1,...,r-1\rbrace$. \\

\noindent
\textbf{Construction of $G_1\diamond(H_0,H_1,...,H_r)$: }\\
The graph is constructed with the vertex set and edge set as follows:\\

$V(G_1\diamond(H_0,H_1,...,H_r))=V(G_1)\cup V(H_0)\cup$ $V(H_1)\cup...\cup V(H_r)$\\

\begin{figure}[!htb]
			\centering
			\includegraphics[width=.5\linewidth]{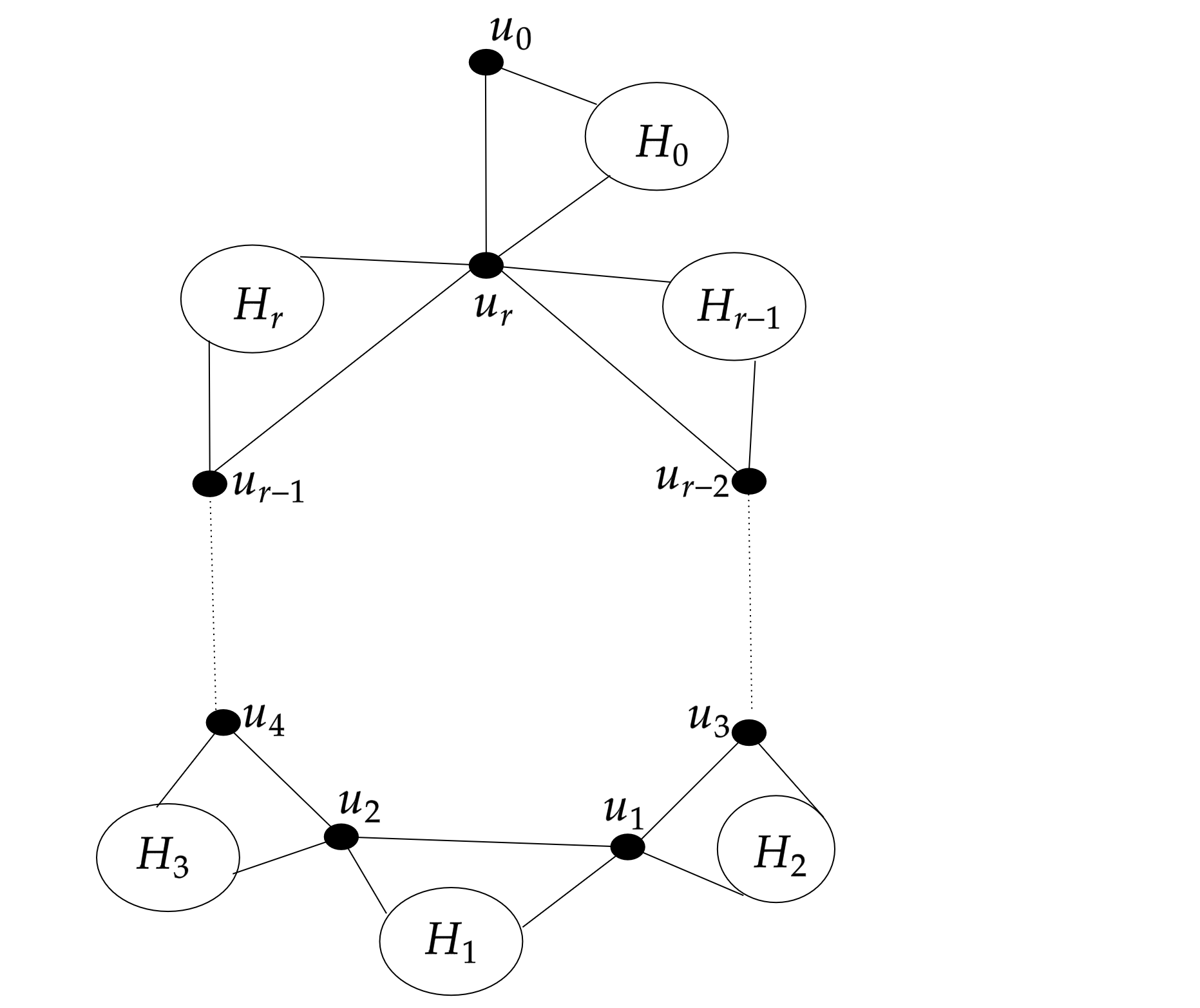}
			\caption{A general representation of $G_1\diamond (H_0,H_1,...,H_r)$}\label{fig:3}
   \end{figure}

$E(G_1\diamond(H_0,H_1,...,H_r))=E(G_1)\cup  E(H_0) \cup $
$E(H_1)\cup...\cup E(H_r)\cup \lbrace u_0v^0_j, u_rv^0_j:v^0_j\in V(H_0)\rbrace$
$\cup\lbrace u_1v^1_j,u_2v^1_j:v^1_j\in V(H_1)\rbrace \cup \lbrace u_1v^2_j,u_3v^2_j:v^2_j\in$
$V(H_2)\rbrace\cup \lbrace u_2v^3_j,u_4v^3_j:v^3_j\in V(H_3)\rbrace \cup ....\cup$
$\lbrace u_{r-2}v^{r-1}_j,u_rv^{r-1}_j:v^{r-1}_j\in V(H_{r-1})\rbrace\cup\lbrace u_{r-1}v^r_j,u_rv^r_j:v^r_j\in V(H_r)\rbrace$ \\
The Figure \ref{fig:3}  gives a general representation of $G_1\diamond (H_0,H_1,...,H_r)$.
	
\subsection{Main results}\label{sec3}

\begin{figure}
		\centering
		\includegraphics[width=1.05\linewidth, height=0.4\textheight]{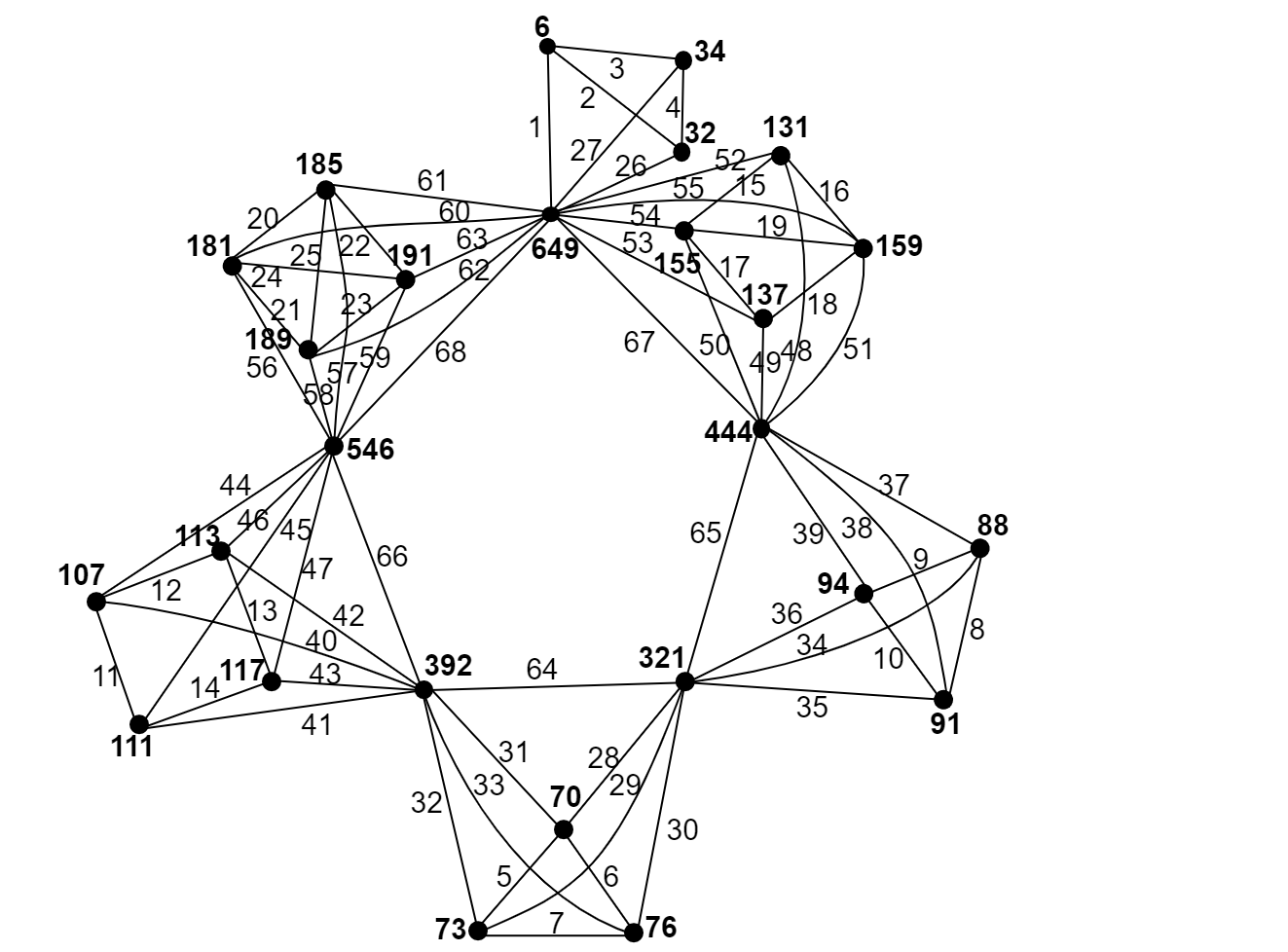}
		\caption{An antimagic labeling of $G_1\diamond (K_2, C_3, C_3, C_4$, a diamond graph, $K_4 )$}
		\label{fig:5}
	\end{figure}

\begin{theorem}\rm
$G_1\diamond (H_0,H_1,...,H_r)$, $r\geq 3$ is antimagic subject to the following conditions:
\begin{equation}
	\left.\begin{aligned}
 &\triangle(H_0)<\delta(H_1);\nonumber\\
 &\triangle(H_i)\leq \delta(H_{i+1}), \hspace{0.1cm}\forall\hspace{0.1cm} i \in \lbrace 1,2,...,r-1\rbrace; \nonumber\\
		&d'(u_0) \leq \delta'(H_i), \hspace{0.1cm} \forall \hspace{0.1cm} i \in \lbrace0,1,...,r\rbrace; \nonumber\\
		& \triangle'(H_r)\leq d'(u_1)\nonumber
	\end{aligned}\right\} (*)
\end{equation}
\begin{proof}\rm 
  We prove that $G_1\diamond (H_0,H_1,...,H_r)$ is antimagic using an algorithmic approach. (Refer Appendix for a better understanding of the following algorithm)
 \begin{algorithm}[!htb]
     \caption{$G_1\diamond (H_0,H_1,...,H_r)$ is antimagic}
  \textbf{Step 1:} Label $E(u_0)$ and $E(H_i)$:\\
   $f(E(u_0))\leftarrow \lbrace 1,2,...,n_0+1 \rbrace$\\
  \For{$i=1,2...,q_0$}{$f(E(H_0))\leftarrow n_0+1+i$}
   \For{$i=1,2,...,q_1$}{$f(E(H_1))\leftarrow n_0+1+q_0+i$}
  \For{$i=1,2,...,q_2$}{$f(E(H_2))\leftarrow n_0+1+q_0+q_1+i$}
   \vdots
  \For{$i=1,2,...,q_r$}{$f(E(H_r))\leftarrow n_0+1+q_0+q_1+...+q_{r-1}+i$}
  Update the PVS: $w'(v^i_j)$, $1\leq j \leq n_i, 0\leq i \leq r$.\\
Update the VS:  $w(u_0)=(1+2+...+n_0)+(n_0+1)$\\
 Rename the vertices $v^i_j$ as $a^i_k$,  $1\leq j,k \leq n_i$, $0\leq i \leq r$ such that $w'(a^i_k)\leq w'(a^i_{k+1})$, $1\leq k \leq n_i-1$, $0\leq i \leq r$.\\
(Observe that $w'(a^i_{j})<w'(a^{i+1}_{k})$, $1\leq j \leq n_i$, $1\leq k \leq n_{i+1}$, $0\leq i \leq r-1$)\\
\textbf{Step 2: } Label the cross edges: \\
Let $c=n_0+1+q_0+q_1+...+q_r$\\
 \For{$i=1,2,...,n_0$}{$f(u_ra^0_i)\leftarrow c+i$}
 \For{$i=1,2,...,n_1$}{$f(u_1a^1_i)\leftarrow c+n_0+i$,\\$f(u_2a^1_i)=c+n_0+n_1+i$}
\end{algorithm}
 
  \begin{algorithm}[!htb]
  \LinesNumbered
	\setcounter{AlgoLine}{21}
	\SetAlgoVlined
   
  \For{$j=1,2,...,r-2$; $i=1,2,...,n_{j+1}$}{$f(u_ja^{j+1}_i) \leftarrow c+n_0+(2n_1+...+2n_j)+i$,\\$f(u_{j+2}a^{j+1}_i)\leftarrow c+n_0+(2n_1+...+2n_j)+n_{j+1}+i$}
  \For{$i=1,2,...,n_r$}{$f(u_{r-1}a^r_i)\leftarrow c+n_0+2n_1+2n_2+...+2n_{r-1}+i$,\\$f(u_ra^r_i)\leftarrow c+n_0+2n_1+2n_2+...+2n_{r-1}+n_r+i$}
Update the VS:\\
$w(a^0_l)=w'(a^0_l)+c+l$, $1\leq l \leq n_0$\\
$w(a^j_i)=w'(a^j_i)+c+n_0+(2n_1+2n_2+...+2n_{j-1})+i+c+n_0+(2n_1+2n_2+...+2n_{j-1})+n_j+i$,  $1\leq j \leq r, 1\leq i \leq n_j$\\
(Observe that $w(a^0_l)< w(a^j_i)< w(a^{j+1}_k)$, $1 \leq l \leq n_0$, $1\leq j \leq r-1, 1\leq i \leq n_j, 1\leq k \leq n_{j+1}$  and $w'(u_i)\leq w'(u_{i+1})$, $1\leq i \leq r-1$)\\
\textbf{Step 3: }Label the remaining edges:\\
Let $b=c+n_0+2n_1+2n_2+...+2n_r$. \\
$f(u_1u_2) \leftarrow b+1$\\
\For{$i=1,2,...,r-2$}{$f(u_iu_{i+2})\leftarrow b+(i+1)$}
$f(u_{r-1}u_r) \leftarrow b+r$\\
Update the VS:\\
$w(u_1)=w'(u_1)+(b+1)+(b+2)$\\
$w(u_j)=w'(u_j)+(b+(j-1))+ (b+(j+1))$, $2\leq j \leq r-1$ \\
$w(u_r)=w'(u_r)+(b+(r-1))+(b+r)$\\
(Observe that $w(u_1)< w(u_j) < w(u_{j+1})<w(u_r)$, $2\leq j \leq r-2$)
\end{algorithm}

\newpage\noindent
\textbf{Proof of distinctness: }In Steps 2 and 3, we have shown the distinctness on the vertex sums $w(a^i_j)$, $1\leq j\leq n_i$, $0\leq i \leq r$ and  $w(u_i)$, $1\leq i \leq r$ respectively. Note that $d'(u_0) \leq \delta'(H_i)$, $ i \in \lbrace 0,1,...,r\rbrace$ (using $(*)$). Clearly, the labels of the edges incident on $u_0$ are less than the labels of the edges incident on $a^i_j$, $1\leq j\leq n_i, 0\leq i \leq r$. Hence,
\begin{align}
  w(u_0)<w(a^i_j), 1\leq j\leq n_i, 0\leq i \leq r \tag{$\alpha$}
\end{align}
 Let\\

\textbf{Set 1:} $\lbrace$$c+n_0+i\hspace{0.1cm} | \hspace{0.1cm} \forall i, \hspace{0.1cm} 1\leq i \leq n_1$$\rbrace$\\
 
	\textbf{Set 2:} $\lbrace$$c+n_0+2n_1+i\hspace{0.1cm} | \hspace{0.1cm} \forall i, \hspace{0.1cm} 1\leq i \leq n_2$$\rbrace$\\
 
	\textbf{Set 3:} $\lbrace$$b+1,b+2$$\rbrace$\\
 
	\textbf{Set 4:} $\lbrace c+n_0+2n_1+2n_2+...+2n_{r-1}+i\hspace{0.1cm} | \hspace{0.1cm} \forall i, \hspace{0.1cm} 1\leq i \leq n_r\rbrace$\\
 
	\textbf{Set 5:} $\lbrace c+n_0+2n_1+2n_2+...+2n_{r-1}+n_r+i \hspace{0.1cm} | \hspace{0.1cm} \forall i, \hspace{0.1cm} 1\leq i \leq n_r\rbrace$\\
 
	\textbf{Set 6:} $\lbrace(n_0+1)+q_0+q_1+...+q_{r-1}+1,...,(n_0+1)+q_0+q_1+...+q_{r-1}+q_r=c\rbrace $\\
 
 Note that $\triangle'(H_r)\leq d'(u_1)$ (using $(*)$). All the vertices in $V(H_r)$ receive $d(v)$ labels of Set 6 where $ v\in V(H_r)$, one of the labels of Set 4, and one of the labels of Set 5 whereas $w(u_1)$ is the sum of the labels of Set 1, Set 2, and Set 3. Hence,

   \begin{align}
      w(a^r_{j}) < w(u_1), 1\leq j \leq n_r \tag{$\beta$}
   \end{align}

Therefore, from $(\alpha),(\beta)$, Step 2, and Step  3, we get,
\begin{align}
 w(u_0)<w(a^i_{j})<w(a^{i+1}_{k})<w(u_1)<...<w(u_r)\nonumber   
\end{align}
  $(1\leq j\leq n_i$, $1\leq k \leq n_{i+1}$ $0\leq i \leq r-1)$.
  Hence, all the vertex sums of the graph $G_1\diamond (H_0,H_1,...,H_r)$ are distinct.
\end{proof}
\end{theorem}
An illustration of the above labeling for the graph $G_1\diamond (H_0,H_1,...,H_5)$ is given in the Figure \ref{fig:5}.
Note that $n_0=2,n_1=3, n_2=3, n_3=4,n_4=4,n_5=4$; $q_0=1,q_1=3,q_2=3,q_3=4,q_4=5, q_5=6$; and the vertex sums  $w(u_0)=6, w(a^0_1)=32, w(a^0_2)=34, w(a^1_1)=70, w(a^1_2)=73, w(a^1_3)=76, w(a^2_1)=88,  w(a^2_2)=91,  w(a^2_3)=94,  w(a^3_1)=107,  w(a^3_2)=111,  w(a^3_3)=113,  w(a^3_4)=117, w(a^4_1)=131, w(a^4_2)=137, w(a^4_3)=155,w(a^4_4)=159, w(a^5_1)=181,w(a^5_2)=185,w(a^5_3)=189,w(a^5_4)=191, w(u_1)=321, w(u_2)=392, w(u_3)=444,w(u_4)=546,w(u_5)=649$.

\subsection{$G_2\diamond (H_1,H_2,...,H_{3p})$ is antimagic}\label{subsec2}
Note that $G_2$ is a spider graph (excluding the spider graph with exactly one vertex of maximum degree 3  containing uneven legs) on $3p+1$ vertices and $3p$ edges containing a vertex $v_0$ of maximum degree $3$ and the other vertices of degree $2$ or $1$ where $p\geq 1$. The graphs $H_1,H_2,...,H_{3p}$ are connected with at least two vertices where $V(H_i)=\lbrace v^i_1,...,v^i_{m_i}\rbrace$ ($|V(H_i)|=m_i$) and $|E(H_i)|=h_i$, $\forall i \in \lbrace 1,2,...,3p\rbrace$. Arrange the graphs $H_i$ satisfying the following condition: $|V(H_i)|\leq |V(H_{i+1})|$, $i \in \lbrace 1,2,..,3p-1\rbrace$. \\

\noindent
\textbf{Construction of $G_2\diamond (H_1,H_2,...,H_{3p})$: }\\
The graph is constructed with the vertex set and edge set as follows:\\

 $V(G_2\diamond (H_1,H_2,...,H_{3p}))=V(G_2)\cup V(H_1)\cup...\cup V(H_{3p})$\\

 $E(G_2\diamond (H_1,H_2,...,H_{3p}))= E(G_2)\cup E(H_1)\cup...\cup E(H_{3p})\cup \mathcal{A} \cup \mathcal{B} \cup \mathcal{C} \cup \mathcal{D}$\\
 
 where $\mathcal{A} = \lbrace x_pv^1_j,x_{p-1}v^1_j:v^1_j\in V(H_1)\rbrace\cup\lbrace x_{p-1}v^4_j,x_{p-2}v^4_j : v^4_j\in V(H_4)\rbrace\cup\lbrace x_{p-2}v^7_j,x_{p-3}v^7_j: v^7_j\in V(H_7)\rbrace\cup...\cup \lbrace x_{1}v^{3p-5}_j,x_{2}v^{3p-5}_j : v^{3p-5}_j\in V(H_{3p-5})\rbrace$\\
 
 $\mathcal{B}= \lbrace y_pv^2_j,y_{p-1}v^2_j:v^2_j\in V(H_2)\rbrace\cup\lbrace y_{p-1}v^5_j,y_{p-2}v^5_j : v^5_j\in V(H_5)\rbrace\cup\lbrace y_{p-2}v^8_j,y_{p-3}v^8_j : v^8_j\in V(H_8)\rbrace\cup...\cup \lbrace y_{1}v^{3p-4}_j,y_{2}v^{3p-4}_j : v^{3p-4}_j\in V(H_{3p-4})\rbrace$\\

$\mathcal{C}=\lbrace z_pv^3_j,z_{p-1}v^3_j:v^3_j\in V(H_3)\rbrace\cup\lbrace z_{p-1}v^6_j,z_{p-2}v^6_j: v^6_j\in V(H_6)\rbrace\cup\lbrace z_{p-2}v^9_j,z_{p-3}v^9_j : v^9_j\in V(H_9)\rbrace\cup...\cup \lbrace z_{1}v^{3p-3}_j,z_{2}v^{3p-3}_j : v^{3p-3}_j\in V(H_{3p-3})\rbrace$\\

$\mathcal{D}=\lbrace x_{1}v^{3p-2}_j, v_0v^{3p-2}_j : v^{3p-2}_j\in V(H_{3p-2})\rbrace \cup\lbrace y_{1}v^{3p-1}_j,v_0v^{3p-1}_j : v^{3p-1}_j\in V(H_{3p-1})\rbrace\cup \lbrace z_{1}v^{3p}_j,v_0v^{3p}_j : v^{3p}_j\in V(H_{3p})\rbrace$

\begin{figure}[!htb]
	\centering
	\includegraphics[width=0.6\textwidth]{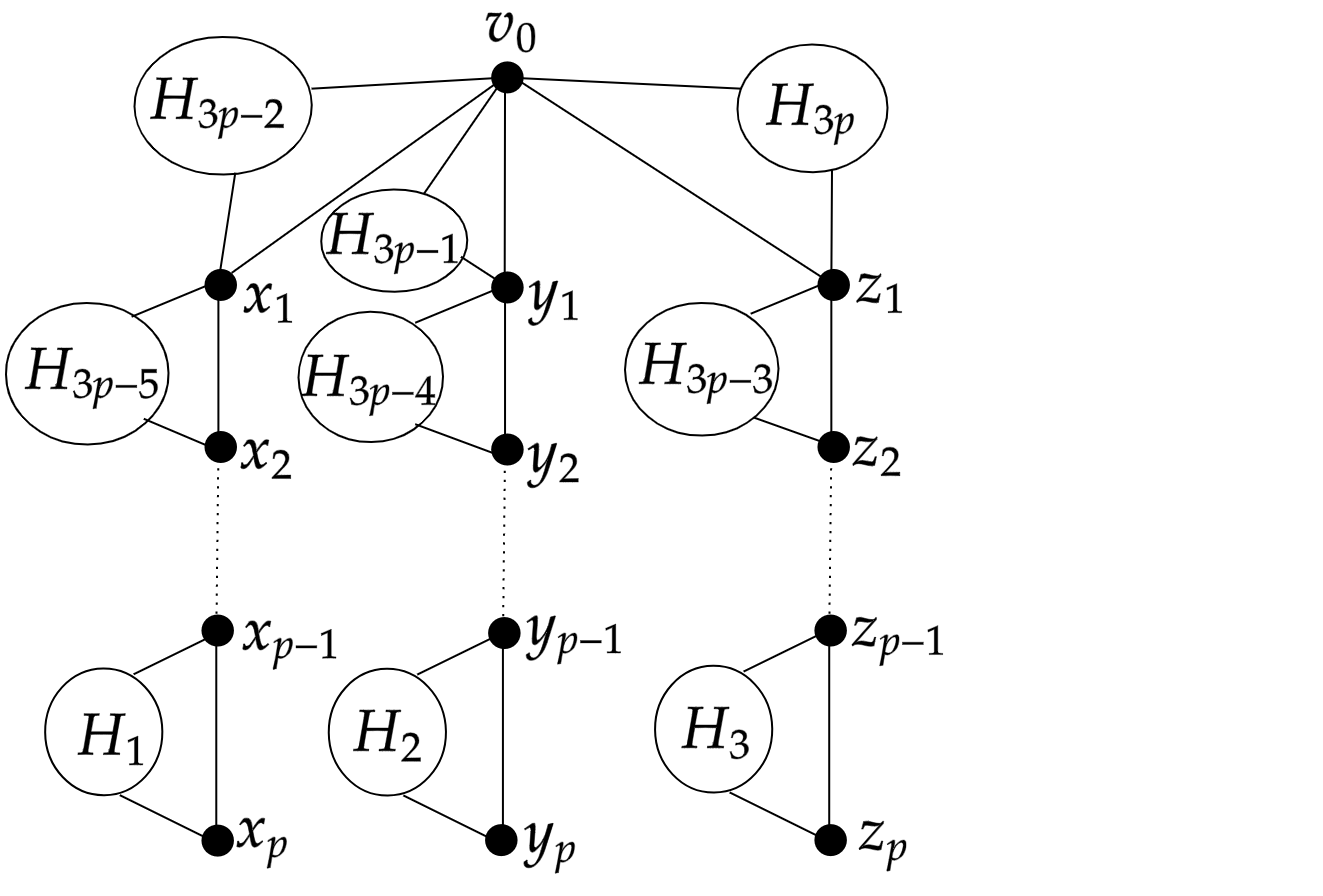} 
	\caption{A general representation of $G_2\diamond (H_1,H_2,...,H_{3p})$}
	\label{fig:6}
\end{figure}

\begin{figure}
		\centering
		\includegraphics[width=1\linewidth, height=0.3\textheight]{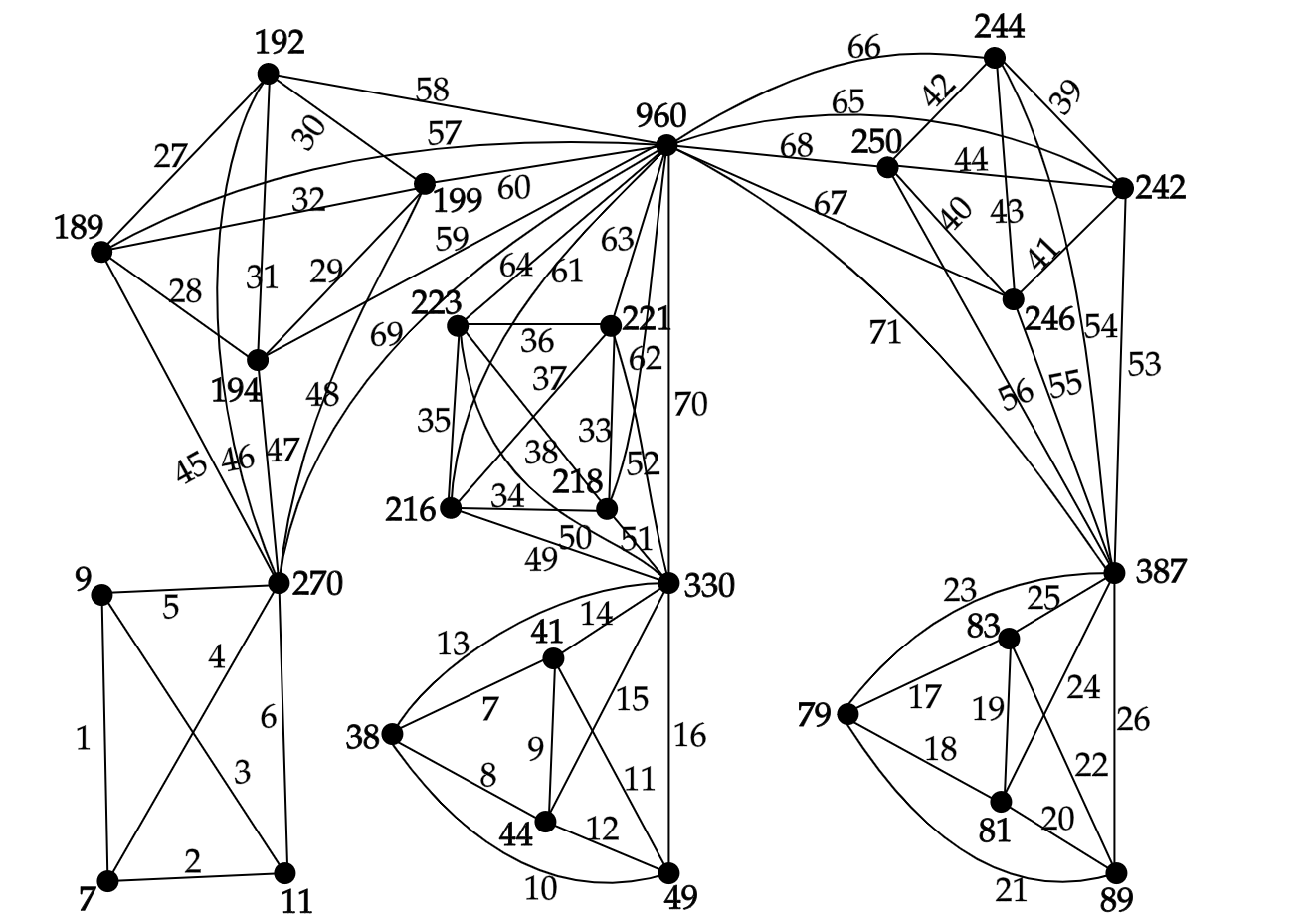}
		 \caption{An antimagic labeling of $G_2\diamond (K_2,C_3,C_3,K_4,K_4,K_4)$}
		\label{fig:8}
	\end{figure}

The general representation of the graph $G_2\diamond (H_1,H_2,...,H_{3p})$ is given in the Figure \ref{fig:6}. Note that the graph $G_2\diamond (H_1,H_2,...,H_{3p})$ has $(m_1+m_2+...+m_{3p})+3p+1$ vertices and $(h_1+h_2+...+h_{3p})+2(m_1+m_2+...+m_{3p})+3p$ edges. The following theorem deals with the graph $G_2$ as a spider containing legs as a path on at most 2 vertices.
\begin{theorem}\rm\label{3.2}
$G_2\diamond (H_1,H_2,...,H_{3p})$, $p=2$ is antimagic subject to the following conditions:
\begin{align}
&\triangle(H_i)\leq \delta(H_{i+1}), \forall\hspace{0.1cm}i \in \lbrace 1,2,...,5\rbrace \nonumber\\ 
& d'(x_{2})\leq \delta'(H_2), \hspace{0.1cm}
d'(y_{2})\leq \delta'(H_3),\hspace{0.1cm}
d'(z_{2})\leq \delta'(H_4)\nonumber
\end{align}
and no restrictions when $p=1$.
\begin{proof}\rm		
To prove the antimagicness of  $G_2\diamond (H_1,H_2,...,H_{3p})$ there are two cases to be discussed: (i) $p=1$ (ii) $p=2$.\\
\noindent
	\textbf{Case(i):} $p=1$\\
 	Here $G_2$ $\cong$ $K_{1,3}$. Since, the maximum degree of the graph $K_{1,3}\diamond (H_1,H_2,H_{3})$ 
 is $|V(K_{1,3}\diamond (H_1,H_2,H_{3}))|-1$,  $K_{1,3}\diamond (H_1,H_2,H_{3})$ is antimagic (refer Lemma 2.1 in \cite{12}). Hence the proof. \\
\noindent
\textbf{Case(ii):} $p=2$\\
 	The labeling technique is represented by an algorithm which is as follows. (Refer Appendix for the better understanding)\\
\begin{algorithm}[!htb]
     \caption{$G_2\diamond (H_1,H_2,...,H_{6})$ is antimagic}
	    \textbf{Step 1:} Label $E(H_1)$ and its cross edges:\\
    \For{$i=1,2,...,h_1$}{ $f(E(H_{1})) \leftarrow i $}
     \For{$i=1,2,...,m_1$}{$f(x_2v^1_i) \leftarrow h_1+i$}
    Update the PVS: $w'(x_{2})$ and $w'(v_i^1)$, $1\leq i \leq m_1$\\
 Rename the vertices $x_{2}$ and $v_i^1$, $1\leq i \leq m_1$ as $a^1_j$,  $1\leq j\leq m_1+1$ such that $w'(a^1_j)\leq w'(a^1_{j+1})$, $1\leq j\leq m_1$.\\
     \For{$j=1,2,...,m_1+1$}{$f(x_{1}a^1_j)=h_1+m_1+j$}
 Update the VS:\\
   $w(a^1_j)=w'(a^1_j)+h_1+m_1+j$, $1\leq j \leq m_1+1$\\
  (Observe that $w(a^1_j)<w(a^1_{j+1})$, $1\leq j \leq m_1$)\\
  \textbf{Step 2: }Label $E(H_2)$ and its cross edges:

\end{algorithm}

 	\begin{algorithm}
   \LinesNumbered
	\setcounter{AlgoLine}{13}
	\SetAlgoVlined
  
   \For{$i=1,2,...,h_2$}{$f(E(H_2)) \leftarrow h_1+2m_1+1+i$}
 \For{$i=1,2,...,m_2$}{$f(y_2v^2_i)\leftarrow h_1+2m_1+1+h_2+i$}
Update the PVS: $w'(y_{2})$ and $w'(v_i^2)$, $1\leq i \leq m_2$. \\
Rename the vertices $y_{2}$ and $v_i^2$, $1\leq i \leq m_2$ as $a^2_j$, $1\leq j \leq m_2+1$ such that $w'(a^2_j)\leq w'(a^2_{j+1})$, $1\leq j\leq m_2$.\\
\For{$j= 1,2,...,m_2+1$}{$f(y_{1}a^2_j)\leftarrow h_1+2m_1+1+h_2+m_2+j$}
Update the VS:\\
$w(a^2_j)=w'(a^2_j)+h_1+2m_1+1+h_2+m_2+j$, $1\leq j \leq m_2+1$\\
(Observe that $w(a^2_j)<w(a^2_{j+1})$, $1\leq j \leq m_2$)\\
\textbf{Step 3: }Label $E(H_3)$ and its cross edges:\\
\For{$i=1,2,..,h_3$}{$f(E(H_3))\leftarrow h_1+2m_1+2+h_2+2m_2+i$}
\For{$i=1,2,...,m_3$}{$f(z_2v^3_i)\leftarrow h_1+2m_1+2+h_2+2m_2+h_3+i$}

Update the PVS: $w'(z_{2})$ and $w'(v_i^3)$, $1\leq i \leq m_3$.\\
Rename the vertices $z_2$ and $v^3_i$, $1\leq i \leq m_3$ as $a^3_j$, $1\leq j \leq m_3+1$ such that $w'(a^3_j)\leq w'(a^3_{j+1})$, $1\leq j\leq m_3$.\\
\For{$j=1,2,...,m_3+1$}{$f(z_{1}a^3_j)=h_1+2m_1+2+h_2+2m_2+h_3+m_3+j$}
Update the VS:\\
$w(a^3_j)=w'(a^3_j)+h_1+2m_1+2+h_2+2m_2+h_3+m_3+j$, $1\leq j \leq m_3+1$\\
(Observe that  $w(a^3_j)<w(a^3_{j+1})$, $1\leq j \leq m_3$)\\
\textbf{Step 4: }  Label $E(H_i)$ and some of its cross edges:\\
Let $z=h_1+2m_1+2+h_2+2m_2+h_3+2m_3+1$.\\
\For{$i=1,2,...,h_4$}{$f(E(H_4))\leftarrow z+i$}
\For{$i=1,2,...,h_5$}{$f(E(H_5))\leftarrow z+h_4+i$}
\For{$i=1,2,...,h_6$}{$f(E(H_6))\leftarrow z+h_4+h_5+i$}
\For{$i=1,2,...,m_4$}{$f(x_1v^4_i)\leftarrow z+h_4+h_5+h_6+i$}
\For{$i=1,2,...,m_5$}{$f(y_1v^5_i)\leftarrow z+h_4+h_5+h_6+m_4+i$}
\For{$i=1,2,...,m_6$}{$f(z_1v^6_i)\leftarrow z+h_4+h_5+h_6+m_4+m_5+i$}
\end{algorithm}

\begin{algorithm}[!htb]
\LinesNumbered
	\setcounter{AlgoLine}{50}
	\SetAlgoVlined
Update the PVS: 
$w'(x_1),w'(y_1),w'(z_1),w'(v_1^{4}),...,w'(v_{m_{4}}^{4}),w'(v_1^{5}),$\\$...,w'(v_{m_{5}}^{5}),w'(v_1^{6}),...,w'(v_{m_{6}}^{6})$.\\
Rename the vertices $x_1,y_1,z_1,v_1^{4},...,v_{m_{4}}^{4},v_1^{5}, ...,v_{m_{5}}^{5},v_1^{6},...,v_{m_{6}}^{6}$ as $c_i$, $1\leq i \leq m_{4}+m_{5}+m_{6}+3$ such that $w'(c_i)\leq w'(c_{i+1})$, $1\leq i \leq m_{4}+m_{5}+m_{6}+2$.\\
 \textbf{Step 5: }Label $E(v_0)$:\\
\For{$i=1,2,...,m_4+m_5+m_6+3$}{$f(v_0c_i)=z+h_4+h_5+h_6+m_{4}+m_{5}+m_{6}+i$}
Update the VS:\\
$w(c_i)=w'(c_i)+z+h_4+h_5+h_6+m_{4}+m_{5}+m_{6}+i$, $1\leq i \leq m_{4}+m_{5}+m_{6}+3$\\
(Observe that $w(c_i)<w(c_{i+1})$, $1\leq i \leq m_{4}+m_{5}+m_{6}+2$)
\end{algorithm}
\noindent
\textbf{Proof of distinctness:}
			Let $A=h_1+2m_1+1$ and $B=A+h_2+2m_2+1$. The condition,
			\begin{align}
				d'(x_{2})\leq \delta'(H_2), &d'(y_{2})\leq \delta'(H_3), d'(z_{2})\leq \delta'(H_4)\nonumber\\
    &\text{ leads to }\nonumber\\
				 d'(a^1_{m_1+1})\leq \delta'(H_2), &d'(a^2_{m_2+1}) \leq \delta'(H_3), 
     d'(a^3_{m_3+1})\leq \delta'(H_4) \text{ respectively.}\nonumber
			\end{align} 
Note that $w(a^1_{m_1+1})=A+\sum\limits_{i=1}^{m_1}(h_1+i)$ and the labels of the edges incident on $a^1_{m_1+1}$ are less than the labels of the edges incident on $a^2_j$, $1\leq j \leq m_2+1$. Observe that $d'(a^2_j)\leq d'(a^2_{m_2+1})$, $1\leq j \leq m_2$.  Hence,
\begin{align}
    w(a^1_{m_1+1})<w(a^2_j), 1\leq j \leq m_2+1
\end{align}

Note that $w(a^2_{m_2+1})=B+\sum\limits_{i=1}^{m_2}(A+h_2+i)$ and the labels of the edges incident on  $a^2_{m_2+1}$ are less than the labels of the edges incident on $a^3_j$, $1\leq j \leq  m_3+1$. Observe that $d'(a^3_{j})\leq d'(a^3_{m_3+1})$, $1\leq j \leq m_3$. Hence,
\begin{align}
    w(a^2_{m_2+1})<w(a^3_j), 1\leq j \leq m_3+1. 
\end{align}

Note that $w(a^3_{m_3+1})=z+\sum\limits_{i=1}^{m_3}(B+h_3+i)$ and the labels of the edges incident on  $a^3_{m_3+1}$ are less than the labels of the edges incident on $c_j$, $1\leq j \leq m_4$. Hence,
\begin{align}
    w(a^3_{m_3+1})<w(c_j), 1\leq j \leq m_4
\end{align}
Therefore, from Steps $1,2,3,5,$ the inequalities $(1), (2),$ and $(3)$, we get 
				\begin{align}
				&w(a^1_j)< w(a^1_{j+1})<w(a^2_k)<w(a^2_{k+1})<w(a^3_l)\nonumber\\
    &<w(a^3_{l+1})<w(c_1)<...<w(c_{m_4+m_5+m_6+3}),\nonumber
			\end{align}
   ($1\leq j\leq m_1$, $1\leq k\leq m_2$, $1\leq l \leq m_3$).
			Also, $w(v_0)=\sum\limits_{i=1}^{m_4+m_5+m_6+3}f(v_0c_i)$. Note that the degree and vertex sum of $v_0$ will be greater than the degree and the vertex sum of any other vertex in the graph $G_2\diamond (H_1,H_2,...,H_{6})$. Hence, $d'(c_{m_4+m_5+m_6+3})<d'(v_0)$ and $w(c_{m_4+m_5+m_6+3})<w(v_0)$. 
			Therefore, all the vertex sums of the graph $G_2\diamond (H_1,H_2,...,H_{6})$ are distinct. Hence proved.
		\end{proof}
  \end{theorem}
An illustration of the above labeling for the graph $G_2\diamond (H_1,H_2,...,H_{6})$ is given in the Figure \ref{fig:8}.
Note that $m_1=2$, $m_2=3$, $m_3=3$, $m_4=4$, $m_5=4$, $m_6=4$, $h_1=1$, $h_2=3$, $h_3=3$, $h_4=6$, $h_5=6$, $h_6=6$, and the vertex sums $w(a^1_1)=7,w(a^1_2)=9, w(a^1_3)=11,w(a^2_1)=38,w(a^2_2)=41,w(a^2_3)=44,w(a^2_4)=49, w(a^3_1)=79,w(a^3_2)=81,w(a^3_3)=83,w(a^3_4)=89, w(c_1)=189,w(c_2)=192,w(c_3)=194,w(c_4)=199,w(c_5)=216,w(c_6)=218,w(c_7)=221,w(c_8)=223,w(c_9)=242,w(c_{10})=244,w(c_{11})=246,w(c_{12})=250,w(c_{13})=270,w(c_{14})=330,w(c_{15})=387,w(v_0)=960$.
\begin{theorem}\rm
$G_2\diamond (H_1,H_2,...,H_{3p})$, $p>2$ is antimagic subject to the following conditions: 
\begin{align}
			&\text{(i) }\triangle(H_i)\leq \delta(H_{i+1}),\forall\hspace{0.1cm}i \in \lbrace1,2,...,3p-1 \rbrace \nonumber\\
			& \text{(ii) }	d'(x_{p})\leq \delta'(H_2), \hspace{0.1cm}
			d'(y_{p})\leq \delta'(H_3), \hspace{0.1cm}
			d'(z_{p})\leq \delta'(H_4)\nonumber\\
			& \text{(iii) }\triangle'(H_i)\leq |V(H_4)|+1, \text{ where $H_i$ is adjacent}\nonumber\\
			& \hspace{0.8cm}\text{to $z_1$ and $z_2$} \nonumber\\
			& \text{(iv) }d'(z_2)\leq \delta'(H_i), \text{ where $H_i$ is adjacent}\nonumber\\
			& \hspace{0.8cm}\text{ to $x_1$ and $ v_0$ }\nonumber
		\end{align}
\begin{proof}
	We follow the Algorithm 2 till Step 3 by replacing $x_2,y_2,z_2$ with $x_p,y_p,z_p$ respectively and $x_1,y_1,z_1$ with $x_{p-1},y_{p-1},z_{p-1}$ respectively. (Refer Appendix  for the better understanding) 
 \begin{algorithm}[!htb]
 \caption{$G_2\diamond (H_1,H_2,...,H_{3p})$ is antimagic, $p>2$ }
Repeat the Steps 1 to 3 from Algorithm 2.\\
\textbf{Step 4: }Label $E(H_i)$:\\
\For{$i=1,2,...,h_4$}{$f(E(H_4))\leftarrow z+i$}
\For{$i=1,2,...,h_5$}{$f(E(H_5))\leftarrow z+h_4+i$}
\For{$i=1,2,...,h_6$}{$f(E(H_6))\leftarrow z+h_4+h_5+i$}
\vdots
\For{$i=1,2,...,h_{3p-3}$}{$f(E(H_{3p-3}))\leftarrow z+h_4+h_5+...+h_{3p-4}+i$}
\textbf{Step 5: }Label the cross edges:\\
Let $L= z+h_4+h_5+...+ h_{3p-3}$. \\
\end{algorithm}

\begin{algorithm}
 \LinesNumbered
	\setcounter{AlgoLine}{12}
	\SetAlgoVlined
 
\For{$i=1,2,...,m_4$}{$f(x_{p-2}v^4_i)\leftarrow L+i$}
\For{$i=1,2,...,m_7$}{$f(x_{p-3}v^7_i)\leftarrow L+m_4+m_5+m_6+i$}
\vdots
\For{$i=1,2,...,m_{3p-5}$}{$f(x_{1}v^{3p-5}_i)\leftarrow L+(m_4+...+m_{3p-6})+i$;}
\For{$i=1,2,...,m_5$}{$f(y_{p-2}v^5_i)\leftarrow L+m_4+i$}
\For{$i=1,2,...,m_8$}{$f(y_{p-3}v^8_i)\leftarrow L+(m_4+...+m_7)+i$}
\vdots
\For{$i=1,2,...,m_{3p-4}$}{$f(y_{1}v^{3p-4}_i)\leftarrow L+(m_4+...+m_{3p-5})+i$;}
\For{$i=1,2,...,m_6$}{$f(z_{p-2}v^6_i)\leftarrow L+m_4+m_5+i$}
\For{$i=1,2,...,m_9$}{$f(z_{p-3}v^9_i)\leftarrow L+(m_4+...+m_8)+i$}
\vdots
\For{$i=1,2,...,m_{3p-3}$}{$f(z_{1}v^{3p-3}_i)\leftarrow L+(m_4+...+m_{3p-4})+i$}
Update the PVS: $w'(v^i_j)$, $4\leq i \leq 3p-3$, $1\leq j \leq m_i$.\\
Rename the vertices $v^i_j$ as $a^i_k$, $1\leq j,k \leq m_i$, $4\leq i \leq 3p-3$ such that
$w'(a^i_k)\leq w'(a^i_{k+1})$, $1\leq k \leq m_i-1, 4\leq i \leq 3p-3$\\
(exclude $w'(x_{p-2}),w'(x_{p-3}),...,w'(x_1)$; $w'(y_{p-2}),w'(y_{p-3}),...,w'(y_1)$; $w'(z_{p-2}),w'(z_{p-3}),...,w'(z_1)$)\\
\textbf{Step 6: } Label the cross edges:\\
Let $N= L+m_4+m_5+...+m_{3p-3}$. \\
\For{$i=1,2,...,m_4$}{$f(x_{p-1}a^4_i)\leftarrow N+i$}
\For{$i=1,2,...,m_7$}{$f(x_{p-2}a^7_i)\leftarrow N+m_4+m_5+m_6+i$}
\vdots
\For{$i=1,2,...,m_{3p-5}$}{$f(x_2a^{3p-5}_i)\leftarrow N+(m_4+...+m_{3p-6})+i$;}

\end{algorithm}

\begin{algorithm}
 \LinesNumbered
	\setcounter{AlgoLine}{41}
	\SetAlgoVlined
 \For{$i=1,2,...,m_5$}{$f(y_{p-1}a^5_i)\leftarrow N+m_4+i$}
 \For{$i=1,2,...,m_8$}{$f(y_{p-2}a^8_i)\leftarrow N+(m_4+...+m_7)+i$}
\vdots
\For{$i=1,2,...,m_{3p-4}$}{$f(y_2a^{3p-4}_i)\leftarrow N+(m_4+...+m_{3p-5})+i$;}
\For{$i=1,2,...,m_6$}{$f(z_{p-1}a^6_i)\leftarrow N+m_4+m_5+i$}

\For{$i=1,2,...,m_9$}{$f(z_{p-2}a^9_i)\leftarrow N+(m_4+...+m_8)+i$}
\vdots
\For{$i=1,2,...,m_{3p-3}$}{$f(z_2a^{3p-3}_i)\leftarrow N+(m_4+...+m_{3p-4})+i$}
Update the VS:\\
$w(a^{i}_j)=w'(a^{i}_j)+N+(m_4+...+m_{i-1})+j$, $4\leq i \leq 3p-3$, $1\leq j \leq m_{i}$\\
(Observe that $w(a^{i}_j)<w(a^{i}_{j+1})$, $4\leq i \leq 3p-3, 1\leq j \leq m_{i}-1$ and $w(a^4_{j_1})<w(a^5_{j_2})<....<w(a^{3p-3}_{j_s})$, $j_1\in \lbrace 1,2,...,m_4\rbrace$, $j_2\in \lbrace 1,2,...,m_5\rbrace$,...., $j_s\in\lbrace 1,2,...,m_{3p-3}\rbrace$)\\
Update the PVS: $w'(x_i),w'(y_i), w'(z_i)$, $2\leq i \leq p-1$ (exclude $w'(x_1),w'(y_1),w'(z_1)$).\\
\textbf{Step 7: } Label the edges present in the legs $s_i$:\\
Let $S= N+m_4+...+m_{3p-3}$.\\
$f(x_{p-1}x_{p-2})\leftarrow S+1$\\
$f(x_{p-2}x_{p-3})\leftarrow S+4$\\
\vdots
$f(x_1x_2)\leftarrow S+(3p-8)$;\\
$f(y_{p-1}y_{p-2})\leftarrow S+2$\\
$f(y_{p-2}y_{p-3})\leftarrow S+5$\\
\vdots
$f(y_1y_2)\leftarrow S+(3p-7)$;\\
$f(z_{p-1}z_{p-2})\leftarrow S+3$\\
$f(z_{p-2}z_{p-3})\leftarrow S+6$\\
\vdots
$f(z_1z_2)\leftarrow S+(3p-6)$;\\
Update the VS:\\

\end{algorithm}

\begin{algorithm}
\LinesNumbered
\setcounter{AlgoLine}{69}
\SetAlgoVlined

\uIf{ $p\neq 3$}{$w(x_{p-1})=w'(x_{p-1})+f(x_{p-1}x_{p-2})$\\ 
$w(x_{p-2})=w'(x_{p-2})+ f(x_{p-1}x_{p-2})+f(x_{p-2}x_{p-3})$\\
 $w(x_{p-3})=w'(x_{p-3})+f(x_{p-2}x_{p-3})+f(x_{p-3}x_{p-4})$\\
 $\vdots$\\
 $w(x_2)=w'(x_2)+f(x_2x_3)+f(x_2x_1)$  }
  \Else{
   $w(x_2)=w'(x_{2})+f(x_{2}x_{1})$ 
  }
  
 (Observe that $w(x_{p-1})<w(x_{p-2})<...<w(x_2)$)\\
\uIf{$p\neq 3$}{
$w(y_{p-1})=w'(y_{p-1})+f(y_{p-1}y_{p-2})$\\
$w(y_{p-2})=w'(y_{p-2})+f(y_{p-1}y_{p-2})+f(y_{p-2}y_{p-3})$\\
$w(y_{p-3})=w'(y_{p-3})+f(y_{p-2}y_{p-3})+f(y_{p-3}y_{p-4})$\\
$\vdots$\\
$w(y_2)=w'(y_2)+f(y_2y_3)+f(y_2y_1)$
  }
\Else{
   $w(y_2)=w'(y_{2})+f(y_{2}y_{1})$ 
  }
  
(Observe that $w(y_{p-1})<w(y_{p-2})<...<w(y_2)$)\\
   \uIf{$p\neq 3$}{
   $w(z_{p-1})=w'(z_{p-1})+f(z_{p-1}z_{p-2})$\\
   $w(z_{p-2})=w'(z_{p-2})+f(z_{p-1}z_{p-2})+f(z_{p-2}z_{p-3})$\\
   $w(z_{p-3})=w'(z_{p-3})+f(z_{p-2}z_{p-3})+f(z_{p-3}z_{p-4})$\\
  $ \vdots$\\
   $w(z_2)=w'(z_2)+f(z_2z_3)+f(z_2z_1)$
    }
     \Else{
  $w(z_2)=w'(z_{2})+f(z_{2}z_{1})$ \
  } 
  
   (Observe that $w(z_{p-1})<w(z_{p-2})<...<w(z_2)$ and $w(x_{p-1})<w(y_{p-1})<w(z_{p-1})<w(x_{p-2})<w(y_{p-2})<w(z_{p-2})<...<w(x_2)<w(y_2)<w(z_2)$)\\
  \textbf{Step 8: }Label $E(H_i)$ and some of its cross edges:\\
  Let $X=S+(3p-6)+h_{3p-2}+h_{3p-1}+h_{3p}$.\\
  Let $M=m_{3p-2}+m_{3p-1}+m_{3p}+3$. \\
   \For{$i=1,2,...,h_{3p-2}$}{$f(E(H_{3p-2}))\leftarrow S+(3p-6)+i$}
  \For{$i=1,2,...,h_{3p-1}$}{$f(E(H_{3p-1}))\leftarrow S+(3p-6)+h_{3p-2}+i$} 
 \end{algorithm}
\newpage
 \begin{algorithm}[!htb]
  \LinesNumbered
	\setcounter{AlgoLine}{103}
	\SetAlgoVlined
 \For{$i=1,2,...,h_{3p}$}{$f(E(H_{3p}))\leftarrow S+(3p-6)+h_{3p-2}+h_{3p-1}+i$}
  \For{$i=1,2,...,m_{3p-2}$}{$f(x_1v^{3p-2}_i)\leftarrow X+i$}
  \For{$i=1,2,...,m_{3p-1}$}{$f(y_1v^{3p-1}_i)\leftarrow X+m_{3p-2}+i$}
 \For{$i=1,2,...,m_{3p}$}{$f(z_1v^{3p}_i)\leftarrow X+m_{3p-2}+m_{3p-1}+i$}
 Update the PVS: $w'(x_1),w'(y_1),w'(z_1),w'(v_1^{3p-2}),...,w'(v_{m_{3p-2}}^{3p-2}),$\\$w'(v_1^{3p-1}),...,w'(v_{m_{3p-1}}^{3p-1}),w'(v_1^{3p}),...,w'(v_{m_{3p}}^{3p})$.\\
 Rename the vertices $x_1,y_1,z_1,v_1^{3p-2},...,v_{m_{3p-2}}^{3p-2}, v_1^{3p-1},...,v_{m_{3p-1}}^{3p-1},$\\$v_1^{3p},...,v_{m_{3p}}^{3p}$ as $c_i$, $1\leq i \leq M$ such that $w'(c_i)\leq w'(c_{i+1})$, $1\leq i \leq M-1$.\\
    \textbf{Step 9: }Label $E(v_0)$:\\
    \For{$i=1,2,...,M$}{$f(v_0c_i)\leftarrow X+m_{3p-2}+m_{3p-1}+m_{3p}+i$}
Update the VS:
$w(c_i)=w'(c_i)+X+m_{3p-2}+m_{3p-1}+m_{3p}+i$, $1\leq i \leq M$\\
(Observe that $ w(c_i)< w(c_{i+1})$, $1\leq i \leq M-1$)
\end{algorithm}

\noindent
\textbf{Proof of distinctness:}
The condition (ii),
\begin{align}
				d'(x_{p})\leq \delta'(H_2), & d'(y_{p})\leq \delta'(H_3), d'(z_{p})\leq \delta'(H_4)\nonumber\\
    &\text{ leads to }\nonumber\\
				d'(a^1_{m_1+1})\leq \delta'(H_2),& d'(a^2_{m_2+1}) \leq \delta'(H_3),
    d'(a^3_{m_3+1})\leq \delta'(H_4) \text{  respectively.}\nonumber
			\end{align}
 We follow the same way as in the proof of distinctness of Theorem \ref{3.2} to obtain the inequalities (1) and (2). \\
   
   Note that $w(a^3_{m_3+1})=z+\sum\limits_{i=1}^{m_3}(B+h_3+i)$ and the labels of the edges incident on  $a^3_{m_3+1}$ are less than the labels of the edges incident on $a^4_j$, $1\leq j \leq m_4$. Hence,
\begin{align}
w(a^3_{m_3+1})<w(a^4_j), 1\leq j \leq m_4
\end{align} 
The condition (iii), $\triangle'(H_i)\leq |V(H_4)|+1$ where $H_i$ is adjacent to $z_1$ and $z_2$, leads to $\triangle'(H_{3p-3})\leq m_4+1$.\\

Observe that $d'(x_{p-1})>m_4+1$ (since $x_{p-1}$ is adjacent to $V(H_4)$ and $V(H_1)$). Hence, $\triangle'(H_{3p-3})<d'(x_{p-1})$. Let \\

	\textbf{Set A: }$\lbrace z+h_4+h_5+...+h_{3p-4}+1,...,	 L\rbrace$\\
 
	\textbf{Set B: }$\lbrace L+m_4+m_5+...+m_{3p-4}+1,...,N\rbrace$\\
 
	\textbf{Set C: }$\lbrace N+m_4+...+m_{3p-4}+1,..., S \rbrace$\\

\noindent
All the vertices in $V(H_{3p-3})$ receive $d(v)$ labels of Set A where $v\in V(H_{3p-3})$, one of the labels of Set B, and one of the labels of Set C whereas $w(x_{p-1}) = (S+1)+\sum\limits_{i=1}^{m_1+1}(h_1+m_1+i)+\sum\limits_{i=1}^{m_4}(N+i)$. Hence, 
\begin{align}
    w(a^{3p-3}_{j}) < w(x_{p-1}), 1\leq j \leq m_{3p-3}
\end{align}
The condition (iv), $d'(z_2)\leq \delta'(H_i)$ where $H_i$ is adjacent to $x_1$ and $v_0$ leads to $d'(z_2)\leq$ $\delta'(H_{3p-2})$.\\

Clearly, the labels of the edges incident on $z_2$ are less than the labels of the edges incident on  $c_j$, $1\leq j \leq m_{3p-2}$. Hence,
\begin{align}
    w(z_2)<w(c_j), 1\leq j \leq m_{3p-2}
\end{align}
Therefore, from Steps $1,2,3,6,7,9,$ the inequalities $(1),(2),(4), (5)$, and $(6)$, we get
\begin{align}
&w(a^1_j)<w(a^1_{j+1})<w(a^2_k)<w(a^2_{k+1})<w(a^3_l)<\nonumber\\
&w(a^3_{l+1})<w(a^4_{j_1})<w(a^5_{j_2})<....<w(a^{3p-3}_{j_s})<\nonumber\\
&w(x_{p-1})<w(y_{p-1})<w(z_{p-1})<....<w(x_2)<\nonumber\\
&w(y_2)<w(z_2)<w(c_1)< ... <w(c_{M})\nonumber
\end{align}
$(1\leq j\leq m_1$, $1\leq k\leq m_2$, $1\leq l \leq m_3$, $1\leq j_1\leq m_4$, $1\leq j_2 \leq m_5$,...,$1\leq j_s \leq m_{3p-3}$).
Also, $w(v_0)=\sum\limits_{i=1}^{M}f(v_0c_i)$. Note that the degree and the vertex sum of $v_0$ will be greater than the degree and the vertex sum of any other vertex in the graph $G_2\diamond (H_1,H_2,...,H_{3p})$. Hence, $d'(c_M)<d'(v_0)$ and $w(c_M)<w(v_0)$. Therefore, all the vertex sums of the graph $G_2\diamond (H_1,H_2,...,H_{3p})$ are distinct. 
\end{proof}
\end{theorem}
  An illustration of the above labeling for the graph $G_2\diamond (H_1,H_2,...,H_{12})$ is given in the Figure \ref{fig:7}.
\begin{figure}
	\centering
	\includegraphics[width=1\linewidth, height=0.40\textheight]{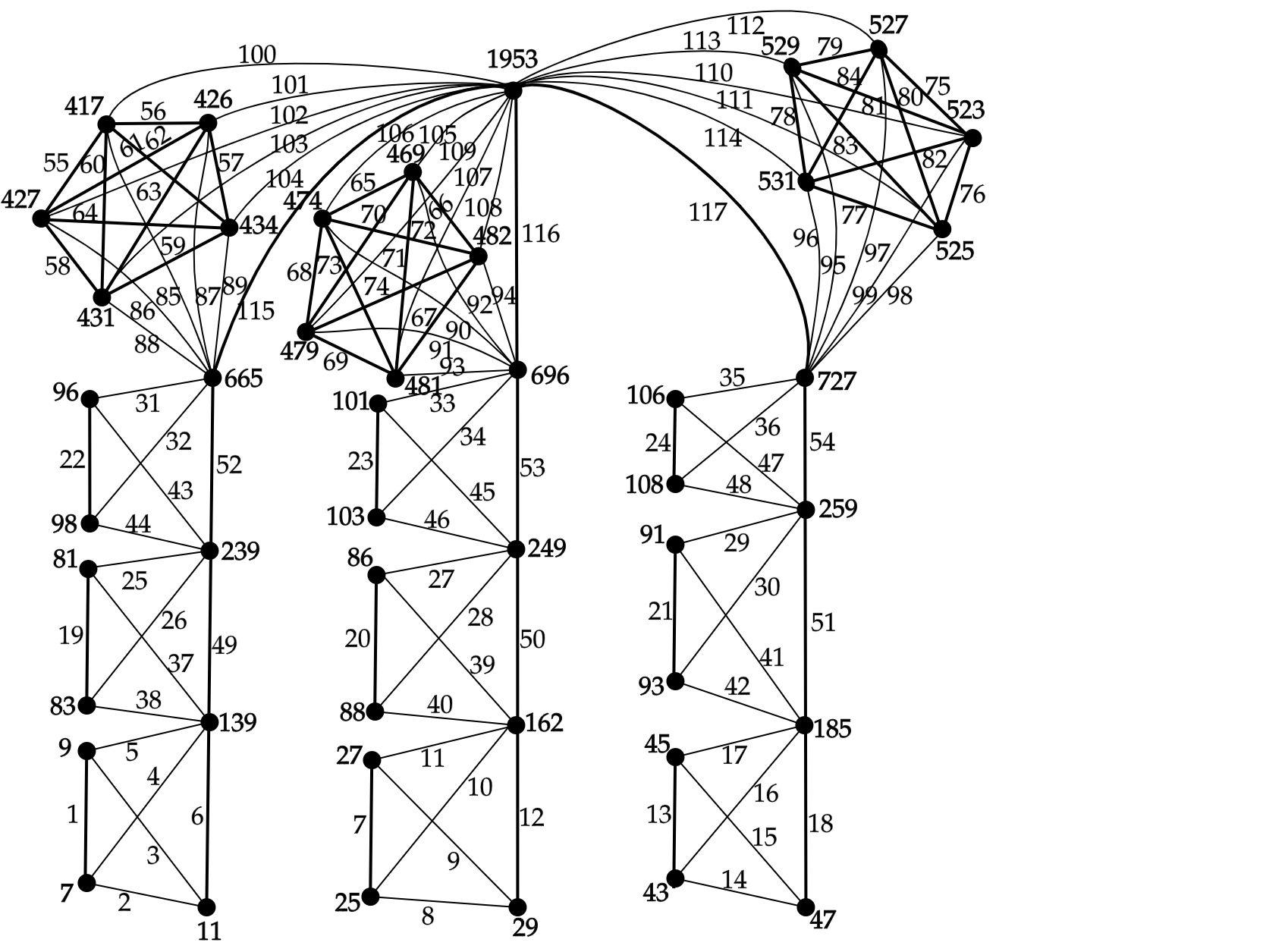}
	\caption{An antimagic labeling of $G_2\diamond (K_2,K_2,K_2,K_2,K_2,K_2,K_2,K_2,K_2,K_5,K_5,K_5)$, ($p=4$)}
	\label{fig:7}
\end{figure}
Note that $m_i=2$, $i\in \lbrace 1,2,...,9\rbrace$, $m_j=5$, $j\in \lbrace 10,11,12\rbrace$, $h_i=1$, $i\in \lbrace 1,2,...,9\rbrace$, $h_j=10$, $j\in \lbrace 10,11,12\rbrace$ and the vertex sums $w(a^1_1)=7,w(a^1_2)=9, w(a^1_3)=11,w(a^2_1)=25,w(a^2_2)=27,w(a^2_3)=29, w(a^3_1)=43,w(a^3_2)=45,w(a^3_3)=47, w(a^4_1)=81,w(a^4_2)=83,w(a^5_1)=86,w(a^5_2)=88,w(a^6_1)=91,w(a^6_2)=93,w(a^7_1)=96,w(a^7_2)=98,w(a^8_1)=101,w(a^8_2)=103,w(a^9_1)=106,w(a^9_2)=108,w(x_2)=139,w(y_2)=162,w(z_2)=185,w(x_3)=239,w(y_3)=249,w(z_3)=259,w(c_1)=417,w(c_{2})=426,w(c_3)=427,w(c_4)=431,w(c_5)=434,w(c_6)=469,w(c_7)=474,w(c_8)=479,w(c_9)=481,w(c_{10})=482,w(c_{11})=523,w(c_{12})=525,w(c_{13})=527,w(c_{14})=529,w(c_{15})=531,w(c_{16})=665,w(c_{17})=696,w(c_{18})=727,w(v_0)=1953$. 

\section{Conclusion}\label{sec13}
Though most of the research focus was on antimagic labeling of general graphs and various products of graphs, there have been no results on antimagic labeling of generalized edge corona of graphs until now. Hence, we were concerned with antimagic labeling of $G\diamond (H_1,H_2,...,H_{m})$  where $|E(G)|=m$ under certain restrictions. In our future work, we plan to focus on the antimagicness of the spider graphs with the maximum degree $3$ having uneven legs. In addition, we pose the following problem as a future direction of research: \lq For any connected graph $G$ with exactly one vertex of maximum degree three, is $G\diamond(H_1,H_2,...,H_{m})$ antimagic?\rq
\\

\end{document}